\documentclass[12pt,leqno]{article}
\usepackage{indentfirst}

\usepackage{amsmath}
\usepackage{amscd}
\usepackage{amssymb}
\usepackage{amsthm}

\def\zR{\ensuremath{\mathbb{R}}}

\def\zN{\ensuremath{\mathbb{N}}}

\numberwithin{equation}{section}
\renewcommand{\theequation}{\thesection.\arabic{equation}}
\makeatletter
\newtheoremstyle{estiloobs}
  {6pt}
  {0pt}
  {\upshape}
  {}
  {}
  {}
  {.5em}
  {(\thmnumber{#2}) \thmname{\textbf{#1}}
   \thmnote{({#3})}}
\makeatother
\theoremstyle{estiloobs}
\newtheorem{remark}[equation]{Remark:}

\makeatletter
\newtheoremstyle{estilotheorem}
  {6pt}
  {0pt}
  {\slshape}
  {}
  {}
  {}
  {.5em}
  {(\thmnumber{#2}) \thmname{\textbf{#1}}
   \thmnote{({#3})}}
\makeatother
\theoremstyle{estilotheorem}
\newtheorem{definition}[equation]{Definition:}

\newtheorem{lemma}[equation]{Lemma:}
\newtheorem{corollary}[equation]{Corollary:}
\newtheorem{theorem}[equation]{Theorem:}
\makeatletter
\newtheoremstyle{estilolista}
  {0pt}
  {0pt}
  {\slshape}
  {}
  {}
  {}
  {.5em}
  {(\thmnumber{#2})\thmname{\textbf{#1}}
   \thmnote{({#3})}}
\makeatother
\theoremstyle{estilolista}

\newcommand{\pedro}{\theequation}

\begin{document}

\title{Weighted inequalities and pointwise estimates for the multilinear fractional integral and
maximal operators}

\vskip 0.3 truecm

\author{Gladis Pradolini \thanks{Research
supported by Consejo Nacional de Investigaciones Cient\'\i ficas y
T\'ecnicas de la Rep\'ublica Argentina and Universidad Nacional
del Litoral.\newline \indent Keywords and phrases: Multilinear
operators, fractional integrals, maximal operators, weighted norm
inequalities.
\newline \indent 1991 Mathematics Subject Classification: Primary
42B25.\newline }}

\date{\vspace{-1.5cm}}

\maketitle

\begin{abstract}
In this article we prove weighted norm inequalities and pointwise
estimates between the multilinear fractional integral operator and
the multilinear fractional maximal. As a consequence of these
estimations we obtain weighted weak and strong inequalities for
the multilinear fractional integral operator. In particular, we
extend some results given in \cite{CPSS} to the multilinear
context. On the other hand we prove weighted pointwise estimates
between the multilinear fractional maximal operator ${\cal
M}_{\alpha,B}$ associated to a Young function $B$ and the
multilinear maximal operators ${\cal M}_{\psi}={\cal M}_{0,\psi}$,
$\psi(t)=B(t^{1-\alpha/(nm)})^{{nm}/{(nm-\alpha)}}$. As an
application of these estimate we obtain a direct proof of the
 $L^p-L^q$ boundedness results of ${\cal M}_{\alpha,B}$ for the
case $B(t)=t$ and $B_k(t)=t(1+\log^+t)^k$ when $1/q=1/p-\alpha/n$.
We also give sufficient conditions on the weights involved in the
boundedness results of ${\cal M}_{\alpha,B}$ that generalizes
those given in \cite{M} for $B(t)=t$. Finally, we prove some
boundedness results in Banach function spaces for a generalized
version of the multilinear fractional maximal operator.
\end{abstract}

\section{Introduction and preliminaries}\label{intro}

An important problem in Analysis is to control certain integral
type operators by means of adequate maximal operators. This
control is sometimes understood in the norm of the spaces where
these operators act. For example, an interesting result due to
Coifman (\cite{C}) establishes that, if $T$ is a
Calder\'on-Zygmund integral operator, $M$ is the Hardy-Littlewood
maximal function and $0<p<\infty$, then the inequality
\begin{equation*}
\int_{\zR^n}|T(f)(x)|^p\, dx \leq C\int_{\zR^n} |Mf(x)|^p\, dx
\end{equation*}
holds for some positive constant $C$. Thus, the maximal function
$M$ controls the singular integral in $L^p$-norm and the
boundedness properties of $M$ in $L^p$-spaces give the boundedness
properties of $T$. The weighted version for $A_{\infty}$ weights
of inequality above is also true (see \cite{CF}).

For the fractional integral operator $I_{\alpha}$, $0<\alpha<n$,
$w\in A_{\infty}$ and $0<p<\infty$, Muckenhoupt and Wheeden
(\cite{MW}) proved the following control-type inequalities
involving the fractional maximal operator $M_{\alpha}$
\begin{equation*}
\int_{\zR^n}|I_{\alpha}(f)(x)|^p w(x)\, dx \leq C\int_{\zR^n}
|M_{\alpha}f(x)|^p w(x)\, dx,
\end{equation*}
and
\begin{equation*}
\sup_{\lambda>0}\lambda^q w(\{I_{\alpha}f>\lambda\})\leq C
\sup_{\lambda>0}\lambda^q w(\{M_{\alpha}f>\lambda\}),
\end{equation*}
where $C$ depends on the $A_{\infty}$-constant of $w$. Then, by
the weighted boundedness results of $M_{\alpha}$, they obtained
the corresponding weighted boundedness results for $I_{\alpha}$.

Similar problems for other operators such that conmutators of
singular and fractional integral operators, non linear
commutators, potential operators, multilinear Calder\'on-Zygmund
operators and multilinear fractional integrals, were studied by
several authors (see for example \cite{P4}, \cite{P5}, \cite{PT},
\cite{CUF}, \cite{BHP}, \cite{CPSS}, \cite{LOPTT} and \cite{M}).
Particularly, in \cite{CPSS}, the authors obtain the boundedness
of the fractional integral operator in term of the fractional
maximal operator in weighted weak $L^1$-spaces, and then, by the
weighted weak boundedness of $M_{\alpha}$, they obtain weighted
weak estimates for $I_{\alpha}$.

Related to the control of commutators of singular and fractional
integral operators appear the iterations of the Hardy- Littlewood
maximal operator $M$ and the composition of the fractional maximal
operator with iterations of $M$. These types of maximal operators
were proved to be equivalent to certain maximal operators
associated to a given Young function (see, for example \cite{P4},
\cite{P5}, \cite{CUF} and \cite{BHP}). Then, the study of the
boundedness properties of these particular maximal operators seem
to be an important tool because they inclose information about the
behaviour of the commutators that they control.

In the multilinear context, there were an increasing interest in
investigate how to control integral operators by maximal
functions. In \cite{GT} the authors proved that the multilinear
Calder\'on-Zygmund operators are controled in $L^p$-norms by the
product of $m$ Hardy-Littlewood maximal operators and they asked
themselves if this product is optimal in some sense. This problem
is then solved in \cite{LOPTT}, where the authors give a strictly
smaller maximal operator and develop a corresponding weighted
theory.

Later, in \cite{M}, a complete study of the weighted boundedness
properties for the multilinear fractional integral operator is
given, and the author proved that this operator is bounded in norm
by the corresponding version of the fractional maximal operator
which generalizes the maximal operator given in \cite{LOPTT}.
Again, the boundedness properties of the ``maximal controller"
gives the boundedness properties of the``controlled operator".

Pointwise estimates between operators are also of interest because
they allow us to obtain boundedness properties of a given operator
by means of the properties of others. For example, related to the
fractional maximal operator and the Hardy-Littlewood maximal
operator a pointwise estimate is given in \cite{CCUF}. Other known
pointwise estimates between the fractional integral operator and
maximal operators are due to Welland and Hedberg (see \cite{W} and
\cite{H}).

In this paper we give ``control type results" for the multilinear
fractional maximal and integral operators. These results involved
pointwise estimates and norm estimates between these operators, of
the type described above. In particular, we extend some results
given in \cite{CPSS} to the multilinear context. On the other hand
we introduce the multilinear fractional maximal operator ${\cal
M}_{\alpha,B}$ associated to a Young function $B$ and we prove
weighted pointwise estimates between these operators and the
multilinear maximal operators ${\cal M}_{\psi}={\cal M}_{0,\psi}$,
where $\psi$ is a given Young function that depends on $B$. As an
application of these estimates we obtain a direct proof of the
 $L^p-L^q$ boundedness results of ${\cal M}_{\alpha,B}$ for the
case $B(t)=t$ and $B_k(t)=t(1+\log^+t)^k$ when $1/q=1/p-\alpha/n$.
We also give sufficient conditions on the weights involved in the
boundedness results of ${\cal M}_{\alpha,B}$ that generalizes
those given in \cite{M} for $B(t)=t$. The importance of a weighted
theory for this maximal function is due to the fact that this
operators are in intimate relation with the commutators of
multisublinear fractional integral operators, as we shall see in a
next paper.

 On the other hand, we study boundedness results in
Banach function spaces for a generalized version of the
multilinear fractional maximal operator involved certain
essentially nondecreasing function $\varphi$.

The paper is organized as follows. In section $\S\ref{section2}$
we statement the main results of this article. We also include
some corollaries and different proofs of results proved in
\cite{M}. The proof of the main results are in $\S\ref{section4}$.
In $\S\ref{section3}$ we give some auxiliary lemmas and finally,
in $\S\ref{section5}$ we define a generalized version of the
multilinear fractional maximal operator and we give some
boundedness estimates in the setting of Banach function spaces.

\bigskip

Before stating the main results of this article, we give some
standard notation. Throughout this paper  $Q$ will denote a cube
in $\zR^n$ with sides parallel to the coordinate axes. With ${\cal
D}$ we will denote the family of dyadic cubes in $\zR^n$.

By a weight we understand a non negative measurable function.

We say that a weight $w$ satisfies a Reverse H\"{o}lder's
inequality with exponent $s$, $RH(s)$, if there exists a positive
constant $C$ such that
\begin{equation*}
\left(\frac{1}{|Q|}\int_Q w^s\right)^{1/s}\leq \frac{w(Q)}{|Q|}.
\end{equation*}

\medskip

By $RH_{\infty}$ we mean the class of weights $w$ such that the
inequality
\begin{equation*}
\sup_{x\in Q}w(x)\leq \frac{C}{|Q|}\int_Q w,
\end{equation*}
holds for every $Q\subset \zR^n$ and some positive constant $C$.
It is easy to check that $RH_{\infty}\subset A_{\infty}$.

\medskip

Now we summarize a few facts about Orlicz spaces. For more
information see \cite{KR} or \cite{RR}.

We say that $B:[0,\infty)\rightarrow [0,\infty)$ is a Young
function if there exists a nontrivial, non-negative and increasing
function $b$ such that $B(t)=\int_0^t b(s)\, ds$. Then $B$ is
continuous, convex, increasing and satisfies $B(0)=0$ and
$\lim_{t\rightarrow\infty}B(t)=\infty$. Moreover, it follows that
$B(t)/t$ is increasing.

Let $B:[0,\infty)\rightarrow [0,\infty)$ be a Young function. The
Orlicz space $L_{B}=L_{B}(\zR^n)$ consists of all measurable
functions $f$ such that for some $\lambda>0$,
\begin{equation*}
\int_{\zR^n}B(|f|/\lambda)<\infty.
\end{equation*}
The space $L_{B}$ is a Banach space endowed with the Luxemburg
norm
\begin{equation*}
\|f\|_{B}=\|f\|_{L_B}=\inf\{\lambda>0:\int_{\zR^n}
B(|f|/\lambda)<\infty\}.
\end{equation*}
The $B$-average of a function $f$ over a cube $Q$ is defined by
\begin{equation*}
\|f\|_{B,Q}=\inf\{\lambda>0:\frac{1}{|Q|}\int_Q B(|f|/\lambda)\leq
1\}.
\end{equation*}
When $B(t)=t$, $\|f\|_{B,Q}=\frac{1}{|Q|}\int_Q |f|$.

We shall say that $B$ is doubling if there exists a positive
constant $C$ such that $B(2t)\leq CB(t)$ for every $t\ge0$. Each
Young function $B$ has an associated complementary Young function
$\tilde{B}$ satisfying
\begin{equation*}
t\leq B^{-1}(t)\tilde{B}^{-1}(t)\leq 2t,
\end{equation*}
for all $t>0$. There is a generalization of H\"{o}lder's
inequality
\begin{equation}\label{holder}
\frac{1}{|Q|}\int_Q |fg|\leq \|f\|_{B,Q}\|g\|_{\tilde{B},Q}.
\end{equation}

A further generalization of H\"{o}lder's inequality (see \cite{O})
is the following: If $\cal A$, $\cal B$ and $\cal C$ are Young
functions such that
\begin{equation*}
{\cal A}^{-1}(t){\cal B}^{-1}(t)\leq{\cal C}^{-1}(t),
\end{equation*}
then
\begin{equation*}
\|fg\|_{{\cal C},Q}\leq 2\|f\|_{{\cal A},Q}\|g\|_{{\cal B},Q}.
\end{equation*}

\begin{definition} Let $0<\alpha<nm$ and $\vec{f}=(f_1,\dots, f_m)$. We
define the multilineal fractional maximal operator
$\mathcal{M}_{\alpha, B}$ associated to a Young function $B$ by
\begin{equation}
\mathcal{M}_{\alpha,B}{\vec{f}}(x)= \sup_{Q\ni
x}|Q|^{\alpha/n}\prod_{i=1}^{m}\|f_i\|_{B,Q}
\end{equation}
where the supremum is taken over all cubes $Q$ containing $x$.
\end{definition}
Even though ${\cal M}_{\alpha,B}$ is sublinear in each entry, we
shall refer to it as the multilineal fractional maximal operator.

\medskip

For $\alpha=0$ we denote $\mathcal{M}_{0, B}=\mathcal{M}_{B}$.
When $B(t)=t$, $\mathcal{M}_{\alpha}=\mathcal{M}_{\alpha,B}$ is
the multilinear fractional maximal operator defined in \cite{M} by
\begin{equation}
\mathcal{M}_{\alpha}{\vec{f}}(x)= \sup_{Q\ni
x}|Q|^{\alpha/n}\prod_{i=1}^{m}\frac{1}{|Q|}\int_Q |f_i|.
\end{equation}
$\mathcal{M}_{0}=\mathcal{M}$ is the multilinear maximal operator
defined in \cite{LOPTT}. When $m=1$ we write $M$ and $M_{\alpha}$
to denote the Hardy-Littlewood and the fractional maximal
operators defined, for a locally integrable function $f$ and
$0<\alpha<n$, by

\begin{equation}\label{hl}
Mf(x)=\sup_{Q\in x}\frac{1}{|Q|}\int_{Q}|f(y)|\, dy
\end{equation}
and
\begin{equation}\label{mf}
M_{\alpha}f(x)=\sup_{Q\in
x}\frac{1}{|Q|^{1-\alpha/n}}\int_{Q}|f(y)|\, dy
\end{equation}
respectively.

\medskip

If we take $g\equiv 1$ in inequality $(\ref{holder})$ it follows
that for every Young function $B$, every $\alpha$ such that
$0\leq\alpha<nm$, the inequality
\begin{equation*}
{\cal M}_{\alpha}(\vec{f})(x)\leq {\cal M}_{\alpha,B}(\vec{f})(x)
\end{equation*}
holds for every $x\in \zR^n$.

\bigskip

The following class of weights was introduced in \cite{LOPTT} and
is a generalization of the Muckenhoupt $A_p$ classes, $p>1$. We
use the notation $\vec{P}=(p_1,\dots,p_m)$.
\begin{definition}
Let $1\leq p_i<\infty$ for $i=1,\dots,m$,
$\frac{1}{p}=\sum_{i=1}^{m}\frac{1}{p_i}$. For each $i=1,\dots ,m$
let $w_i$ be a weight and $\vec{w}=(w_1,\dots ,w_m)$. We say that
$\vec{w}$ satisfies the $A_{\vec{P}}$ condition if
\begin{equation}\label{apmulti}
\sup_Q
\left(\frac{1}{|Q|}\int_Q\left(\prod_{i=1}^{m}w_i^{p/p_i}\right)\right)^{1/p}\prod_{i=1}^{m}
\left(\frac{1}{|Q|}\int_Q w_i^{1-p_i'}\right)^{1/p_i'}<\infty.
\end{equation}
When $p_i=1$, $\left(\frac{1}{|Q|}\int_Q
w_i^{1-p_i'}\right)^{1/p_i'}$ is understood as $(\inf_Q
w_i)^{-1}$.
\end{definition}

\bigskip

Condition $(\ref{apmulti})$ is called the multilinear
$A_{\vec{P}}$ condition.

\bigskip

\section{Statement of the main
results}\label{section2}

In this section we establish the main results of this article. For
a sake of completeness we consider subsections.

\bigskip

\centerline {\sl Pointwise estimates for ${\cal M}_{\alpha,B}$}

\bigskip

For $0<\alpha<nm$ let $B$ be a Young function such that
$t^{\frac{\alpha}{nm}}B^{-1}(t^{1-\frac{\alpha}{nm}})\leq
B^{-1}(t)$. Let $\psi$ be the function defined by
$\psi(t)=B(t^{1-\alpha/(nm)})^{{nm}/{(nm-\alpha)}}$. From lemma
$\ref{Younggamma}$ bellow, $\psi$ is a Young function. The
following result gives a pointwise estimate between the
multilinear fractional maximal associated to the Young function
$B$, $\mathcal{M}_{\alpha, B}$ and the multilineal maximal
operator $\mathcal{M}_{\psi}$ associated to the Young function
$\psi$, and is an useful tool to obtain boundedness results for
$\mathcal{M}_{\alpha, B}$.

\bigskip

\begin{lemma}\label{malphapeso} Let $0<\alpha<nm$. Let $B$ be a Young function such
that
\begin{equation}\label{condi1}
t^{\frac{\alpha}{nm}}B^{-1}(t^{1-\frac{\alpha}{nm}})\leq B^{-1}(t)
\end{equation}
and $\psi(t)=B(t^{1-\alpha/(nm)})^{{nm}/{(nm-\alpha)}}$.

For each $i=1,\dots,m$, let $p_i$, $q_i$ and $s_i$ be the real
numbers defined, respectively, by $1\leq p_i<nm/\alpha$,
$\frac{1}{q_i}=\frac{1}{p_i}-\frac{\alpha}{nm}$ and
$s_i=\left(1-\alpha/(nm)\right)q_i$  and $p$, $q$ and $s$ be the
real numbers given by $\frac{1}{p}=\sum_{i=1}^{m}\frac{1}{p_i}$,
$\frac{1}{q}=\sum_{i=1}^{m}\frac{1}{q_i}$ and
$\frac{1}{s}=\sum_{i=1}^{m}\frac{1}{s_i}$. Let $w_1,\dots,w_m$ be
$m$ weights. If $\vec{f_w}=(f_1/w_1,\dots, f_m/w_m)$ and
$\vec{g}=(f_1^{p_1/s_1}w_1^{-q_1/s_1},\dots,
f_m^{p_m/s_m}w_m^{-q_m/s_m})$  then
\begin{equation}\label{puntualpesada}
\mathcal{M}_{\alpha, B}\vec{f_w}(x)\leq
\mathcal{M}_{\psi}{\vec{g}}(x)^{1-\alpha/(nm)}\left(\prod_{i=1}^{m}\|f_i\|_{L^{p_i}}^{p_i}\right)^{\frac{\alpha
}{nm}}.
\end{equation}

\end{lemma}

\bigskip

\begin{remark}When $B(t)=t$ we have that
$\psi(t)=t$. Then, from inequality $(\ref{puntualpesada})$, we get
the following pointwise estimate between the multilinear
fractional maximal operator $\mathcal{M}_{\alpha}$ and the
multilineal maximal operator $\mathcal{M}$

\begin{equation}\label{puntualpesada1} \mathcal{M}_{\alpha}\vec{f_w}(x)\leq
\mathcal{M}{\vec{g}}(x)^{1-\alpha/(nm)}\left(\prod_{i=1}^{m}\|f_i\|_{L^{p_i}}^{p_i}\right)^{\frac{\alpha
}{nm}}.
\end{equation}
In the case $m=1$ the result above was obtained in \cite{GPS}.
\end{remark}

\bigskip

\begin{remark}For $0<\alpha<nm$ and $k\in \zN$ let $B_k$ be the
Young function defined by $B_k(t)=t(1+\log^+ t)^k$. Then $B_k$
satisfies $(\ref{condi1})$. Let
$\psi_k(t)=B_k(t^{1-\alpha/(nm)})^{{nm}/{(nm-\alpha)}}\cong
t(1+\log^+ t)^{k{nm}/{(nm-\alpha)}}$. From the lemma above we get
the following pointwise estimate
\begin{equation}\label{puntualpesada2} \mathcal{M}_{\alpha,
L(\log L)^k}\vec{f_w}(x)\leq \mathcal{M}_{L(\log
L)^{\frac{knm}{nm-\alpha}}}{\vec{g}}(x)^{1-\alpha/(nm)}\left(\prod_{i=1}^{m}\|f_i\|_{L^{p_i}}^{p_i}\right)^{\frac{\alpha
}{nm}}.
\end{equation}
\end{remark}

\vspace{0.5cm}

\centerline{\sl Weighted boundedness results for ${\cal
M}_{\alpha,B}$}

\bigskip

 As an easy consequence of inequality $(\ref{puntualpesada1})$ and
 the weighted boundedness results for the multilinear maximal operator
 $\mathcal{M}$ proved in \cite{LOPTT} we obtain a direct proof of
 the weighted weak and strong boundedness of the multilinear fractional maximal operator
 $\mathcal{M}_{\alpha}$ proved in
 \cite{M}, when $p$ and $q$ satisfy $1/q=1/p-\alpha/n$ and $1<p_i<nm/\alpha$, $i=1,\dots,m$. Actually, in \cite{M} the author proves that the conditions
 on the weights are also necessary (see theorems 2.7 and 3.6 in \cite{M} applied to this
 case). These results are given in the following two theorems.

\bigskip

\begin{theorem}\label{debil}
Let $0<\alpha<nm$ and let $p_i$ and $q$ be defined as in lemma
\ref{malphapeso}. Let $\vec{f}=(f_1,\dots, f_m)$, If $(u,\vec{w})$
satisfy
\begin{equation}\label{weak}
\sup_Q \left(\frac{1}{|Q|}\int_Q
u\right)^{1/q}\prod_{i=1}^{m}\left(\frac{1}{|Q|}\int_Q
w_i^{-pi'}\right)^{1/p_i'}<\infty
\end{equation}
then
\begin{equation*}
\|\mathcal{M}_{\alpha}\vec{f}\|_{L^{q,\infty}(u)}\leq
C\prod_{i=1}^{m}\|f_i w_i\|_{L^{p_i}}.
\end{equation*}
\end{theorem}

\bigskip

\begin{theorem}\label{lplq}
Let $0<\alpha<nm$ and let $1<p_i<nm/\alpha$ and $q_i$, $s_i$, $p$,
$q$, and $s$ be defined as in lemma \ref{malphapeso}. Let
$\vec{w^{q}}=(w_1^{q_1},\dots,w_m^{q_m})$. If $\vec{f}=(f_1,\dots,
f_m)$, $\vec{S}=(s_1,\dots, s_m)$ and $\vec{w^q}\in A_{\vec{S}}$,
then
\begin{equation*}
\|\mathcal{M}_{\alpha}\vec{f}\, (\Pi_{i=1}^{m}w_i)\|_{L^q}\leq C
\prod_{i=1}^{m}\|f_i w_i\|_{L^{p_i}}.
\end{equation*}
\end{theorem}

\begin{remark}
It is easy to check that $\vec{w^{q}}\in A_{\vec{S}}$ if and only
if $\vec{w}=(w_1,\dots,w_m)$ belongs to the $A_{\vec{P},q}$
classes introduced in \cite{M}. This equivalence is a
generalization to the multilinear case of that proved by
Muckenhoupt and Wheeden in the linear case, which establishes that
a weight $w \in A_{p,q}$ if and only if $w^q \in A_s$ with $1\leq
p<n/\alpha$, $1/q=1/p-\alpha/n$ and $s=1+q/p'$. For more details
see \cite{MW}.
\end{remark}

\vspace{0.8cm}

The following corollary is a consequence of theorem $\ref{debil}$
applied to the weights $u=\prod_{i=1}^m u_i^{q/q_i}$ and $w_i=
M(u_i)^{1/q_i}$, where $M$ is the Hardy-Littlewood maximal
operator.

\bigskip

\begin{corollary}\label{debil2}
Let $0<\alpha<nm$ and let $p_i$, $q_i$, $s_i$, $p$, $q$, and $s$
be defined as in lemma \ref{malphapeso}. Let $\vec{f}=(f_1,\dots,
f_m)$ and $u=\prod_{i=1}^{m}u_i^{q/q_i}$ then
\begin{equation*}
\|\mathcal{M}_{\alpha}\vec{f}\|_{L^{q,\infty}(u)}\leq
C\prod_{i=1}^{m}\|f_i M(u_i)^{1/q_i}\|_{L^{p_i}}.
\end{equation*}
\end{corollary}

\bigskip

From the weak and strong characterizations obtained in [M,
Theorems 2.7 and 2.8] applied to the case $p=q$, we obtain the
following result.

\begin{theorem}
Suppose that $0<\alpha<nm$, $1\leq p_1,\dots ,p_m<mn/\alpha$ and
$\frac{1}{p}=\sum_{i=1}^{m}\frac{1}{p_i}$. Let
$u=\prod_{i=1}^{m}u_i^{p/p_i}$ and $v=\prod_{i=1}^{m}v_i^{1/p_i}$.
Then
\begin{equation*}
\|\mathcal{M}_{\alpha}\vec{f}\|_{L^{p,\infty}(u)}\leq
C\prod_{i=1}^{m}\|f_i\|_{L^{p_i}(M_{\alpha p_i/m}(u_i))},
\end{equation*}
and
\begin{equation*}
\|\mathcal{M}_{\alpha}\vec{f}\, v\|_{L^p}\leq
C\prod_{i=1}^{m}\|f_i\,M_{\alpha p_i/m}(v_i)\|_{L^{p_i}},
\end{equation*}
where $M_{\alpha p_i/m}$ denotes the fractional maximal operator
defined in $(\ref{mf})$ with $\alpha$ replaced by $\alpha p_i/m$.
\end{theorem}

\medskip

The proof of the first inequality above follows from the fact that
the weights $u$ and $w_i=M_{\alpha p_i/m}(u_i)$ satisfy the
condition on the weights in [M, Theorem 2.7]. On the other hand,
the weights $v$ and $w_i=M_{\alpha p_i/m}(v_i)^{1/p_i}$ satisfy
the hypotheses in [M, Theorem 2.8] and thus we obtain the second
inequality.

\vspace{0.5cm}

Before state the next result, we introduce the following class of
Young functions related to the boundedness of the sublinear
maximal $M_B$ between Lebesgue spaces. For more information see
\cite{P1}.

\begin{definition}
Let $1<r<\infty$. A Young function $B$ is said to satisfy the
$B_r$ condition if for some constant $c>0$,
\begin{equation*}
\int_c^{\infty}\frac{B(t)}{t^r}\, \frac{dt}{t}<\infty.
\end{equation*}
\end{definition}

\medskip

\begin{theorem}\label{pesosorlicz}
Let $0\leq\alpha<nm$, $1< p_i<\infty$, $i=1,\dots,m$,
$\frac{1}{p}=\sum_{i=1}^{m}\frac{1}{p_i}$. Let $q$ be a real
number such that $1/m<p\leq q<\infty$. Let $B$, $A_i$, and $C_i$,
$i=1,\dots,m$, be Young functions such that
$A_i^{-1}(t)C_i^{-1}(t)\leq B^{-1}(t)$, $t>0$ and $C_i$ is
doubling and satisfies the $B_{p_i}$ condition for every
$i=1,\dots,m$. Let $(\nu,\vec{w})$ weights that satisfy
\begin{equation}\label{orl}
\sup_Q
|Q|^{\alpha/n+1/q-1/p}\left(\frac{1}{|Q|}\int_Q\nu^q\right)^{1/q}
\prod_{i=1}^m \|w_i^{-1}\|_{A_i,Q}<\infty
\end{equation}
then
\begin{equation*}
\|{\cal M}_{\alpha,B}\vec{f}\, \nu\|_{L^q}\leq C
\prod_{i=1}^{m}\|f_i w_i\|_{L^{p_i}}.
\end{equation*}
holds for every $\vec{f}\in L^{p_1}(w_1^{p_1})\times \dots
L^{p_m}(w_m^{p_m})$.
\end{theorem}

\medskip

\begin{remark}
The linear case of theorem above was proved in \cite{CUF}, and in
\cite{CUP} for the case $\alpha=0$ and $p=q$. For $B(t)=t$,
theorem $\ref{pesosorlicz}$ gives two weighted results proved in
\cite{M} for the multilineal fractional maximal operator ${\cal
M}_{\alpha}$. The first one ([M, Theorem 2.8]) is obtained by
considering $A_i=t^{rp_i'}$ and $C_i=t^{(rp_i')'}$ for some $r>1$
and the second ([M, Theorem 2.10]) is obtained by taking
$A_i=t^{p_i'}(1+\log^+t)^{p_i'-1+\delta}$ and
$C_i=\frac{t^{p_i}}{(1+\log^+t)^{1+\delta(p_i-1)}}$ for
$\delta>0$.
\end{remark}

\bigskip

As a consequence of theorem $\ref{pesosorlicz}$ and the pointwise
estimate given in $(\ref{puntualpesada2})$, we obtain the
following result about the boundedness of
$\mathcal{M}_{\alpha,B_k}$ for multilinear weights in the
$A_{\vec{S}}$ class defined above, where $\vec{S}=(s_1,\dots
,s_m)$ and $B_k(t)=t(1+\log^+ t)^k$. In the proof, we also use the
pointwise estimate given in $(\ref{puntualpesada2})$.

\medskip

\begin{corollary}\label{acotacionap}
Let $0\leq \alpha<nm$ and let $p_i$, $p$, $q_i$, $q$, $s_i$ and
$s$ be defined as in lema $\ref{malphapeso}$. For each $k\in \zN$
let $B_k(t)=t(1+\log^+ t)^k$. Let $\vec{w^q}=(w_1^{q_1},\dots
,w_m^{q_m})$. If $\vec{f}=(f_1,\dots, f_m)$ and
$\vec{S}=(s_1,\dots,s_m)$ then the inequality
\begin{equation*}
\|\mathcal{M}_{\alpha,B_k}\vec{f}\, (\Pi_{i=1}^{m}w_i)\|_{L^q}\leq
C\prod_{i=1}^{m}\|f_i w_i\|_{L^{p_i}}
\end{equation*}
holds for every $\vec{f}$ if and only if $\vec{w^q}$ satisfies the
$A_{\vec{S}}$ condition.
\end{corollary}

\bigskip




\bigskip

\centerline{\sl Weighted weak type inequalities for the
multilinear fractional integral operator}

\bigskip

In this section we obtain weighted estimates for the multilinear
fractional maximal and integral operator.

The following definition of the multilinear fractional integral
operator was considered by several authors (see, for example,
\cite{G}, \cite{KS}, \cite{GK} and \cite{M}).

\begin{definition}
Let $0<\alpha<nm$ and $\vec{f}=(f_1,\dots,f_m)$. The multilinear
fractional integral is defined by
\begin{equation*}
{\cal I}_{\alpha}\vec{f}(x)=\int_{(\zR^n)^m}\frac{f_1(y_1)\dots
f_m(y_m)}{(|x-y_1|+\dots+|x-y_m|)^{mn-\alpha}}\, d\vec{y},
\end{equation*}
where the integral in convergent if $\vec{f}\in {\cal S}\times
\dots \times{\cal S}$.
\end{definition}

\bigskip

Particularly, we study weighted weak type inequalities for the
multilinear fractional maximal and integral operator. For the
first one we obtain the following result.

\medskip

\begin{theorem}\label{debilmaximal} Let $0\leq\alpha<nm$, $\vec{w}=(w_1,\dots,w_m)$
and $u=\prod_{i=1}^m w_i^{1/m}$. Then
\begin{equation*}
u(\{x\in \zR^n: {\cal M}_{\alpha}\vec{f}(x)>\lambda^m\})^m\leq
C\prod_{i=1}^m\int_{\zR^n} \frac{|f_i|}{\lambda}\,
M_{\alpha/m}w_i,
\end{equation*}
where $M_{\alpha/m}$ denotes the fractional maximal operator of
order $\alpha/m$ defined in $(\ref{mf})$.
\end{theorem}

\bigskip

The case $\alpha=0$ of theorem above was proved in \cite{LOPTT}.
For $m=1$ this is a well known result proved in \cite{FS}.

\bigskip

In \cite{CPSS} the authors considered the problem of find weights
$W$ such that
\begin{equation*}
w(\{x\in \zR^n : |I_{\alpha}f(x)|>\lambda\})\leq
\frac{C}{\lambda}\int_{\zR^n}|f(x)|W(x)\, dx
\end{equation*}
for a given weight $w$, for every $\lambda>0$ and for suitable
functions $f$. Particularly, they obtain that the weight
$W=M_{\alpha}(M_{L(\log L)^{\delta}}w)$, $\delta>0$, works.
Motivated from the linear case, we study an analogous problem in
the multilinear context and we obtain the following result.

\bigskip

\begin{theorem}\label{controlmaximal}
Let $0<\alpha<nm$, $\delta>0$ and $u=\prod_{i=1}^m w_i^{1/m}$.
Then
\begin{equation}\label{debilialpha}
\|{\cal I}_{\alpha}\vec{f}\|_{L^{1/m,\infty}(u)}\leq
C\prod_{i=1}^{m}\int_{\zR^n}|f_i|M_{\alpha/m}M_{L(\log
L)^{\delta}}(w_i).
\end{equation}
and, in particular
\begin{equation*}
\|{\cal I}_{\alpha}\vec{f}\|_{L^{1/m,\infty}(u)}\leq
C\prod_{i=1}^{m}\int_{\zR^n}|f_i|M_{\alpha/m}M^2(w_i).
\end{equation*}
\end{theorem}

\bigskip

The result above is an immediate consequence of the next theorem.

\bigskip

\begin{theorem}\label{debildebil}
Let $0<\alpha<nm$, $\delta>0$ and let $u$ be a weight. Then

\begin{equation*}
\|{\cal I}_{\alpha}\vec{f}\|_{L^{1/m,\infty}(u)}\leq C\|{\cal
M}_{\alpha}\vec{f}\|_{L^{1/m,\infty}(M_{L(\log L)^{\delta}}(u))}.
\end{equation*}
\end{theorem}

\bigskip

Then, the proof of $(\ref{debilialpha})$ follows by observing that
\begin{equation*}
M_{L(\log L)^{\delta}}(u)=M_{L(\log L)^{\delta}}(\prod_{i=1}^m
w_i^{1/m})\leq \prod_{i=1}^{m}M_{L(\log L)^{\delta}}(w_i)^{1/m},
\end{equation*}
which is a consequence of the generalized H\"{o}lder's inequality
in Orlicz spaces. Then, an application of theorem
$\ref{debilmaximal}$ gives the desired result.

\vspace{0.7cm}

\bigskip

Recall that a weight $v$ satisfies the $RH_{\infty}$ condition if
there exists a positive constant $C$ such that the inequality
\begin{equation*}
\sup_{x\in Q}v(x)\leq \frac{C}{|Q|}\int_Q v
\end{equation*}
holds for every $Q\subset \zR^n$.

\bigskip

\begin{lemma}\label{dospesos}Let $0<\alpha<nm$. Let $v$ be a weight satisfying the $RH_{\infty}$
condition. Then, there exists a positive constant $C$ such that,
if $u=\prod_{i=1}^m w_i^{1/m}$ and $\vec{f}=(f_1,\dots,f_m)$,
\begin{equation*}
\int_{\zR^n}{\cal I}_{\alpha}\vec{f}(x) u(x) v(x) dx \leq
C\int_{\zR^n}{\cal M}_{\alpha}\vec{f}(x)Mu(x) v(x)\,dx,
\end{equation*}
where $M$ is the Hardy-Littlewood maximal function defined in
$(\ref{hl})$.
\end{lemma}

\bigskip

The following theorem establish some kind of control of the
multilinear fractional integral operator by the multilinear
fractional maximal in $L^p$, $0<p\leq 1$.

\bigskip

\begin{theorem}\label{p<1}
Let $0<p\leq 1$ and let $u$ be a weight. Then
\begin{equation*}
\int_{\zR^n}|{\cal I}_{\alpha}\vec{f}(x)|^p u(x) dx \leq
C\int_{\zR^n}|{\cal M}_{\alpha}\vec{f}(x)|^p Mu(x)\,dx.
\end{equation*}
\end{theorem}

\bigskip

In the linear case, lemma $\ref{dospesos}$ and theorem $\ref{p<1}$
were proved in \cite{CPSS}.

\vspace{0.7cm}

\centerline {\sl Pointwise estimates between ${\cal I}_{\alpha}$
and ${\cal M}_{\alpha}$}

\bigskip

A pointwise estimate relating both, the multilinear fractional and
maximal operators is given in the next result.

\medskip

\begin{theorem}{\bf(Welland's type inequality)}\label{Welland} Let
$0<\alpha<nm$ and $0<\epsilon<\min\{\alpha,nm-\alpha\}$. Then, if
$\vec{f}=(f_1,\dots,f_m)$ where $f_i$'s are bounded functions with
compact support, then
\begin{equation*}
|{\cal I}_{\alpha}\vec{f}(x)|\leq
C\left(\mathcal{M}_{\alpha+\epsilon}{\vec{f}(x)}\mathcal{M}_{\alpha-\epsilon}{\vec{f(x)}}\right)^{1/2},
\end{equation*}
where $C$ only depends on $n$, $m$, $\alpha$ and $\epsilon$.
\end{theorem}

\bigskip

In \cite{M}, the author proves the following result.

\medskip

\begin{theorem}{\bf [M, theorem 2.2]}
Suppose that $0<\alpha<nm$, $1<p_1,\dots,p_m<\infty$ and $q$ is a
number that satisfies $1/m<p\leq q<\infty$. Suppose that one of
the two following conditions holds.\\
\noindent {\rm (i)} $q>1$ and $(\nu,\vec{w})$ are weights that
satisfy
\begin{equation*}
\sup_Q
|Q|^{\alpha/n+1/q-1/p}\left(\frac{1}{|Q|}\int_Q\nu^{qr}\right)^{1/(qr)}
\prod_{i=1}^m \left(\frac{1}{|Q|}\int_Q
w_i^{-p_i'r}\right)^{1/(p_i'r)}<\infty
\end{equation*}
for some $r>1$.\\
\noindent {\rm (ii)} $q\leq 1$ and $(\nu,\vec{w})$ are weights
that satisfy
\begin{equation*}
\sup_Q
|Q|^{\alpha/n+1/q-1/p}\left(\frac{1}{|Q|}\int_Q\nu^{q}\right)^{1/q}
\prod_{i=1}^m \left(\frac{1}{|Q|}\int_Q
w_i^{-p_i'r}\right)^{1/(p_i'r)}<\infty
\end{equation*}
for some $r>1$.\\
Then the inequality
\begin{equation*}
\|{\cal I}_{\alpha}\vec{f}\nu\|_q\leq
C\prod_{i=1}^m\|f_iw_i\|_{p_i}.
\end{equation*}
holds for every $\vec{f}\in L^{p_1}(w_1^{p_1})\times\dots
L^{p_m}(w_m^{p_m})$.
\end{theorem}

\bigskip

A direct proof of theorem above for the case $q>1$ can be given
combining theorem $\ref{Welland}$  with theorem
$\ref{pesosorlicz}$ applied to the case $A_i(t)=t^{rp_i'}$, and
proceeding as in the corresponding result in \cite{GCM} (theorem
6.5).

\vspace{0.7cm}

\section{Auxiliary results}\label{section3}

In this section we give some technical lemmas used in the proof of
the main results in this paper.

\begin{lemma}\label{Younggamma}
Let $B$ be a Young function and $0<\gamma<1$. Then
$\psi(t)=B(t^{\gamma})^{1/\gamma}$ is a Young function.
\end{lemma}

\bigskip

\noindent {\it Proof}: It is enough to prove that there exists a
nontrivial, non-negative and increasing function $g$ such that
$\psi(t)=\int_{0}^{t}g(s)\, ds$. This function $g$ is given by
$g(s)=b(s^{\gamma})\left(\frac{B(s^{\gamma})}{s^{\gamma}}\right)^{(1/\gamma)-1}$,
where $b$ is a non-negative and increasing function such that
$B(t)=\int_0^tb(s)ds$. The function $g$ has the desired
properties. $\square$

\vspace{0.7cm}

The next lemma establishes the relation between the dyadic a
non-dyadic multilinear fractional maximal operators. Let ${\cal
M}_{\alpha,B}^k$ be defined as ${\cal M}_{\alpha,B}$ but over
cubes with side length less or equal than $2^k$,
$Q_k=Q(0,2^{k+2})$, $\tau_tg(x)=g(x-t)$ and
$\vec{\tau}_t(\vec{f})=(\tau_tf_1,\dots,\tau_tf_m)$.

\begin{lemma}\label{dyadic}
For each $k$, $\vec{f}$ and every $x\in\zR^n$ and $0<q<\infty$,
there exists a constant $C$, depending only on $n$, $m$, $\alpha$
and $q$ such that
\begin{equation}\label{relation}
{\cal M}_{\alpha,B}^k(\vec{f})(x)^q\leq
\frac{C}{|Q_k|}\int_{Q_k}(\tau_{-t}\circ {\cal
M}_{\alpha,B}^d\circ \vec{\tau}_t)(\vec{f})(x)^q\, dt
\end{equation}
\end{lemma}

For the linear case and $\alpha=0$ this result was proved by
Fefferman and Stein in \cite{FS} and can be also found in
\cite{GCRF}. In the multilinear context and $\alpha=0$ the result
above is given in \cite{LOPTT}, and for $B(t)=t$ and $\alpha>0$,
in \cite{M}. The proof of lemma \ref{dyadic} is an easy
modification of any of the mentioned results and we omit it.

\bigskip

In order to prove theorem $\ref{debildebil}$ we need the following
results. The first of them was proved in \cite{M} for the
multilinear integral operator. For the linear case, a proof can be
found in \cite{P2}.

\bigskip

\begin{lemma}{\bf[M]}
Let $g$ and $f_i$, $i=1,\dots,m$ be positive functions with
compact support and let $u$ be a weight. Then there exists a
family of dyadic cubes $\{Q_{k,j}\}$ and a family of pairwise
disjoint subsets $\{E_{k,j}\}$, $E_{k,j}\subset Q_{k,j}$ with
\begin{equation*}
|Q_{k,j}|\leq C |E_{k,j}|
\end{equation*}
for some positive constant $C$ and for every $k,\, j$ and such
that
\begin{eqnarray}\label{discreta}
\int_{\zR^n}{\cal I}_{\alpha}\vec{f}(x) u(x) g(x) dx&\leq&
C\sum_{k,j}|Q_{k,j}|^{\alpha/n}\left(\frac{1}{|Q_{k,j}|}\int_{Q_{k,j}}u(x)
g(x)\,
dx\right)\\
&&\quad\times
\left(\prod_{i=1}^m\frac{1}{|3Q_{k,j}|}\int_{3Q_{k,j}} f_i(y_i)\,
dy_i\right) |E_{k,j}|.\nonumber
\end{eqnarray}
\end{lemma}

\bigskip

The following lemma was proved in \cite{CN} and gives examples of
weights in the $RH_{\infty}$ class.

\medskip

\begin{lemma}\label{CN}
Let $g$ be any function such that $Mg$ is finite a.e.. Then
$(Mg)^{-\alpha}\in RH_{\infty}$, $\alpha>0$.
\end{lemma}

\bigskip

\section{Proofs}\label{section4}



\noindent {\it Proof of lemma $\ref{malphapeso}$}:

\medskip

The proof is based in some ideas from lemma 2.8 in \cite{GPS}. Let
$g_i$ be a function such that $g_i^{s_i}w_i^{q_i}=f_i^{p_i}$. Then
$f_i/w_i=g_i^{s_i/p_i}w_i^{q_i/p_i-1}=g_i^{s_i/p_i+\alpha/(nm)-1}g_i^{1-\alpha/(nm)}w_i^{(q_i/p_i-1)}$.
Let $r=nm/(nm-\alpha)$ and $r'=nm/\alpha$. If $s$ and $s_i$ are
defined as in the hypotheses of the theorem we get
\begin{equation}\label{uno}
\left(\frac{q_i}{p_i}-1\right)r'=\left(\frac{q_i}{p_i}-1\right)\frac{nm}{\alpha}=q_i
\end{equation}
and
\begin{eqnarray}\label{dos}
\left(\frac{s_i}{p_i}+\frac{\alpha}{nm}-1\right)r'&=&\left(\frac{s_i}{p_i}+
\frac{\alpha}{nm}-1\right)\frac{nm}{\alpha}\\\nonumber
&=&\left(\left(1-\frac{\alpha}{nm}\right)\frac{q_i}{p_i}+\frac{\alpha}{nm}-1\right)\frac{nm}{\alpha}\\\nonumber
&=&\left(1-\frac{\alpha}{nm}\right)\left(\frac{q_i}{p_i}-1\right)\frac{nm}{\alpha}\\\nonumber
&=&\left(1-\frac{\alpha}{nm}\right)q_i\\\nonumber &=& s_i
\end{eqnarray}

Let $B$ and $\psi$ be the functions in the hypotheses of the
theorem. From lemma $\ref{Younggamma}$ $\psi$ is a Young function.
Let $\phi(t)=B(t)^{nm/(nm-\alpha)}$. Then, by the properties of
the function $B$ we obtain
\begin{equation*}
\phi^{-1}(t)\, t^{\alpha/nm}\leq C B^{-1}(t).
\end{equation*}
By applying H\"older's inequality, and using $(\ref{uno})$ and
$(\ref{dos})$ we obtain that
\begin{eqnarray*}
\|f_i/w_i\|_{B,Q}&=& \|g_i^{s_i/p_i}w_i^{q_i/p_i-1}\|_{B,Q}\\
&=& \|g_i^{1-\alpha/nm}g_i^{s_i/p_i+\alpha/nm-1}w_i^{q_i/p_i-1}\|_{B,Q}\\
&\leq& \|g_i^{1-\alpha/nm}\|_{\phi,Q}\|g_i^{s_i/p_i+\alpha/nm-1}w_i^{q_i/p_i-1}\|_{nm/\alpha,Q}\\
&=&\frac{1}{|Q|^{\alpha/nm}}\|g_i\|_{\psi,Q}^{1-\alpha/nm}\|f_i\|_{p_i}^{\alpha
p_i/nm}.
\end{eqnarray*}
where we have used that
$\|g_i^{1-\alpha/nm}\|_{\phi,Q}=\|g_i\|_{\psi,Q}^{1-\alpha/nm}$.
Then

\begin{eqnarray*}
|Q|^{\alpha/n}\prod_{i=1}^{m}\|f_i/w_i\|_{B,Q}
&\leq&\prod_{i=1}^{m}\|g_i\|_{\psi,Q}^{1-\alpha/(nm)}\prod_{i=1}^{m}
\|f_i\|_{p_i}^{\alpha p_i/(nm)}.\\
&\leq&\mathcal{M}_{\psi}\vec{g}(x)^{1-\alpha/(nm)}\left(\prod_{i=1}^{m}
\|f_i\|_{p_i}^{p_i}\right)^{\alpha/(nm)},\\
\end{eqnarray*}
and inequality $(\ref{puntualpesada})$ follows by taking supremum
over the cubes $Q$ in $\zR^n$. $\square$

\bigskip

\noindent {\it Proof of theorem $\ref{debil}$}: We use the same
notation as in the proof of lemma $\ref{malphapeso}$. Thus, it is
enough to prove that
\begin{equation*}
\|\mathcal{M}_{\alpha}\vec{f_w}\|_{L^{q,\infty}(u)}\leq
C\prod_{i=1}^{m}\|f_i\|_{p_i},
\end{equation*}
and then replace $f_i$ by $f_iw_i$.

From the hypotheses on the weights and raising the quantity in
$(\ref{weak})$ to the power $1-\alpha/(nm)$ we obtain that
\begin{equation*}
\sup_Q \left(\frac{1}{|Q|}\int_Q
u\right)^{1/s}\prod_{i=1}^{m}\left(\frac{1}{|Q|}\int_Q
w_i^{-pi'}\right)^{1/s_i'}<\infty
\end{equation*}
o, equivalently
\begin{equation}\label{weak1}
\sup_Q \left(\frac{1}{|Q|}\int_Q
u\right)^{1/s}\prod_{i=1}^{m}\left(\frac{1}{|Q|}\int_Q
w_i^{q_i(1-s_i')}\right)^{1/s_i'}<\infty
\end{equation}

By inequality $(\ref{puntualpesada1})$ and from $(\ref{weak1})$
and the weighted weak boundedness result for ${\cal M}$ proved in
\cite{LOPTT} we obtain that
\begin{eqnarray}\label{condi2}
\|\mathcal{M}_{\alpha}\vec{f_w}\|_{L^{q,\infty}(u)}&\leq&C\left(\prod_{i=1}^{m}\|f_i\|_{p_i}^{p_i}\right)^
{\alpha/nm}
\|\mathcal{M}g\|_{L^{s,\infty}(u)}^{1-\alpha/nm}\\
&\leq&C\left(\prod_{i=1}^{m}\|f_i\|_{p_i}^{p_i}\right)^{\alpha/nm}\left
(\prod_{i=1}^{m}\|g_i\|_{L^{s_i}(w_i^{q_i})}\right)^{1-\alpha/nm}\nonumber\\
&=&C\left(\prod_{i=1}^{m}\|f_i\|_{p_i}^{p_i}\right)^{\alpha/nm}\left
(\prod_{i=1}^{m}\|f_i\|_{p_i}^{p_i/s_i}\right)^{1-\alpha/nm}\nonumber\\
&=&C\prod_{i=1}^{m}\|f_i\|_{p_i},\nonumber
\end{eqnarray}
where we have used that $p_i\,
\alpha/(nm)+(p_i/s_i)(1-\alpha/(nm))=1$. Thus the proof is done.
$\square$

\noindent {\it Proof of theorem $\ref{lplq}$}: Let
$\nu=\prod_{i=1}^m w_i$. As in the proof above , it is enough to
show that
\begin{equation*}
\|\mathcal{M}_{\alpha}\vec{f_w}\nu\|_{q}\leq C
\prod_{i=1}^{m}\|f_i\|_{p_i},
\end{equation*}
but this inequality can be obtained in a similar way to that in
$(\ref{condi2})$ by replacing $\|\mathcal{M}g\|_{L^{s,\infty}(u)}$
by $\|\mathcal{M}g\|_{L^{s}(\nu^q)}$ and then using the weighted
strong boundedness result proved in \cite{LOPTT}. $\square$

\bigskip

\noindent {\it Proof of theorem $\ref{pesosorlicz}$}: We first
consider the dyadic version ${\cal M}_{\alpha,B}^d$ de ${\cal
M}_{\alpha,B}$ defined by
\begin{equation*}
{\cal M}_{\alpha,B}^d=\sup_{Q\in {\cal D}:x\in
Q}|Q|^{\alpha/n}\prod_{i=1}^m \|f_i\|_{B,Q}.
\end{equation*}
where ${\cal D}$ denotes the set of dyadic cubes in $\zR^n$. Let
$a$ be a constant satisfying $a>2^{mn}$ and for each $k$ let
\begin{equation*}
\Omega_k=\{x\in \zR^n:{\cal M}_{\alpha,B}^d(\vec{f})(x)>a^k\}.
\end{equation*}
It is easy to see that an analogue of the Calder\'on Zygmund
decomposition in Orlicz spaces holds for ${\cal M}_{\alpha, B}^d$
and, therefore there is a family of maximal non-overlapping dyadic
cubes $\{Q_{j,k}\}$ such that $\Omega_k=\cup_{j}Q_{j,k}$ and
\begin{equation*}
a^k<|Q_{j,k}|^{\alpha/n}\prod_{i=1}^m\|f_i\|_{B,Q_{k,j}}\leq
2^{nm}a^k.
\end{equation*}
Moreover, each ${\Omega}_{k+1}\subset {\Omega }_k$ and the sets
$E_{k,j}=Q_{k,j}\backslash(Q_{k,j}\cap {\Omega}_{k+1})$ are
disjoint and satisfy
\begin{equation}\label{desi}
|Q_{k,j}|<\beta|E_{k,j}|
\end{equation}
for some $\beta>1$. Then, by the generalized H\"{o}lder's
inequality and condition $(\ref{orl})$ we obtain
\begin{eqnarray*}
\int_{\zR^n}{\cal M}_{\alpha,B}^d(\vec{f})^q \nu^q &=& \sum_k
\int_{\Omega_k
\backslash \Omega_{k+1}}{\cal M}_{\alpha,B}^d(\vec{f})^q\nu^q\\
&\leq&a^q\sum_k a^{kq}\nu^q(\Omega_k)\\
&\leq&a^q\sum_{k,j} a^{kq}\nu^q(Q_{k,j})\\
&\leq&C\sum_{k,j} \left(|Q_{k,j}|^{\alpha/n}\prod_{i=1}^m\|f_i\|_{B,Q_{k,j}}\right)^q\nu^q(Q_{k,j})\\
&\leq&C\sum_{k,j}
\left(|Q_{k,j}|^{\alpha/n}\prod_{i=1}^m\|f_iw_i\|_{C_i,Q_{k,j}}\right)^q\left(\prod_{i=1}^m\|w_i^{-1}
\|_{A_i,Q_{k,j}}^q\right)\nu^q(Q_{k,j})\\
&\leq&C\sum_{k,j}
\left(\prod_{i=1}^m\|f_iw_i\|_{C_i,Q_{k,j}}\right)^q|Q_{k,j}|^{q/p}.\\
\end{eqnarray*}
Now, from the fact that $p\leq q$ and using $(\ref{desi})$, the
multilinear H\"{o}lder's inequality and the hypotheses on $C_i$ we
obtain that
\begin{eqnarray*}
\left(\int_{\zR^n}{\cal
M}_{\alpha,B}^d(\vec{f})^q\nu^q\right)^{1/q}
 &\leq&C\left(\sum_{k,j}
\left(\prod_{i=1}^m\|f_iw_i\|_{C_i,Q_{k,j}}\right)^p|Q_{k,j}|\right)^{1/p}\\
&\leq&C\left(\sum_{k,j}
\left(\prod_{i=1}^m\|f_iw_i\|_{C_i,Q_{k,j}}\right)^p|E_{k,j}|\right)^{1/p}\\
&\leq&C\prod_{i=1}^m\left(\sum_{k,j} \|f_iw_i\|_{C_i,Q_{k,j}}^{p_i}|E_{k,j}|\right)^{1/p_i}\\
&\leq&C\prod_{i=1}^m\left(\int_{\zR^n} M_{C_i}(f_iw_i)^{p_i}\right)^{1/p_i}\\
&\leq&\prod_{i=1}^m \|f_iw_i\|_{L^{p_i}}.
\end{eqnarray*}

To prove the non-dyadic case we use lemma $\ref{dyadic}$. Thus,
from $(\ref{relation})$ it follows that
\begin{equation}\label{normaq}
\|{\cal M}_{\alpha,B}(\vec{f})\,\nu\|_q \leq
\sup_t\|\tau_{-t}\circ{\cal
M}^d_{\alpha,B}\circ\vec{\tau_{t}}(\vec{f})\, \nu\|_q.
\end{equation}

If the weights $(\nu,\vec{w})$ satisfy condition $(\ref{orl})$,
then the weights $(\tau_{t}(\nu),\vec{\tau_t}\vec{w})$ satisfy the
same condition with constant independent of $t$. Then, applying
the dyadic case, we obtain

\begin{eqnarray*}
\|(\tau_{-t}\circ{\cal
M}^d_{\alpha,B}\circ\vec{\tau_{t}})(\vec{f})\, \nu\|_q&=&\|({\cal
M}^d_{\alpha,B}\circ\vec{\tau_{t}})(\vec{f})\, \tau_{t}\nu\|_q\\
&\leq& C \prod_{i=1}^m\|\tau_tf_i\, \tau_tw_i\|_{p_i}\\
&\leq& C \prod_{i=1}^m \|f_iw_i\|_{p_i},\\
\end{eqnarray*}
with $C$ independent of $t$. Then, from $(\ref{normaq})$ we obtain
that
\begin{equation*}
\|{\cal M}_{\alpha,B}(\vec{f})\,\nu\|_q \leq  C \prod_{i=1}^m
\|f_iw_i\|_{p_i}.\square
\end{equation*}

\bigskip

\noindent {\it Proof of corollary $\ref{acotacionap}$}: We begin
by proving the case $\alpha=0$. If $\vec{w^p}\in A_{\vec{P}}$ then
we have that $w_i^{p_i(1-p_i')}=w_i^{-p_i'}\in A_{mp_i'}$ (see
\cite{LOPTT}). Then, for each $i=1,\dots , m$ there exist $s_i>1$
such that $w_i^{-p_i'}$ satisfies a reverse H\"{o}lder inequality
with exponent $s_i$. Let $A_i(t)=t^{s_ip_i'}$ and
$C_{i,k}(t)=\left(t(1+\log^+t)^k\right)^{(s_i p_i')'}$. Then we
have that $A_i^{-1}(t)C_{i,k}^{-1}(t)\cong B_k^{-1}(t)$ and
$C_{i,k}\in B_{p_i}$. Thus, since $\vec{w^p}\in A_{\vec{P}}$ we
obtain
\begin{eqnarray*}
\left(\frac{1}{|Q|}\int_Q(\Pi_{i=1}^{m}w_i)^p\right)^{1/p}\prod_{i=1}^{m}\|w_i^{-1}\|_{A_i,Q}&\leq&
\left(\frac{1}{|Q|}\int_Q(\Pi_{i=1}^{m}w_i)^p\right)^{1/p}\prod_{i=1}^{m}\left(\frac{1}{|Q|}\int_Q
w_i^{-s_ip_i'}\right)^{1/(s_ip_i')}\\
&\leq&\left(\frac{1}{|Q|}\int_Q(\Pi_{i=1}^{m}w_i)^p\right)^{1/p}\prod_{i=1}^{m}\left(\frac{1}{|Q|}\int_Q
w_i^{-p_i'}\right)^{1/(p_i')}\\
&\leq& C.
\end{eqnarray*}
Then by theorem $\ref{pesosorlicz}$ applied to the case $\alpha=0$
and $p=q$ we obtain that
\begin{equation*}
\|{\cal M}_{B_k}(\vec{f})(\Pi_{i=1}^{m}w_i)\|_{L^p}\leq
\prod_{i=1}^{m}\|f_iw_i\|_{L_{p_i}}.
\end{equation*}

The other implication is a consequence of the inequality ${\cal
M}(\vec{f})\leq {\cal M}_{B_k}(\vec{f})$ and the boundedness
results proved in \cite{LOPTT} for the multilinear maximal
operator ${\cal M}$.

 Now we prove the case $\alpha>0$. Let us first suppose that $\vec{w^q}\in A_{\vec{S}}$. It is enough to show
that the following inequality
\begin{equation}\label{general}
\|{\cal M}_{\alpha, B_k}(\vec{f_w})(\Pi_{i=1}^{m}w_i)\|_{L^q}\leq
\prod_{i=1}^{m}\|f_i\|_{L_{p_i}}.
\end{equation}
holds for every $\vec{f}=(f_1,\dots,f_m)$.

In order to prove $(\ref{general})$ we use inequality
$(\ref{puntualpesada2})$. Let
$\psi_k(t)=t(1+\log^+t)^{knm/(nm-\alpha)}$. Then, from the case
$\alpha=0$ we obtain that
\begin{eqnarray*}
\|{\cal M}_{\alpha,
B_k}(\vec{f_w})(\Pi_{i=1}^{m}w_i)\|_{L^q}&\leq&C\|{\cal
M}_{\psi_k}(\vec{g})^{1-\frac{\alpha}{nm}}(\Pi_{i=1}^{m}w_i)\|_{q}
\left(\prod_{i=1}^{m}\|f_i\|_{L^{p_i}}^{\frac{p_i\alpha}{nm}}\right)\\
&\leq& C\|{\cal M}_{L(\log
L)^{\frac{knm}{nm-\alpha}}}(\vec{g})(\Pi_{i=1}^{m}w_i^{q/s})\|_{L^s}^{s/q}
\left(\prod_{i=1}^{m}\|f_i\|_{L^{p_i}}^{\frac{p_i\alpha}{nm}}\right)\\
&\leq& C\|{\cal M}_{L(\log
L)^{\left[\frac{knm}{nm-\alpha}\right]+1}}(\vec{g})(\Pi_{i=1}^{m}w_i^{q/s})\|_{L^s}^{s/q}
\left(\prod_{i=1}^{m}\|f_i\|_{L^{p_i}}^{\frac{p_i\alpha}{nm}}\right)\\
&=& C\|{\cal
M}_{B_{\left[\frac{knm}{nm-\alpha}\right]+1}}(\vec{g})(\Pi_{i=1}^{m}w_i^{q/s})\|_{L^s}^{s/q}
\left(\prod_{i=1}^{m}\|f_i\|_{L^{p_i}}^{\frac{p_i\alpha}{nm}}\right)\\
&\leq&
C\left(\prod_{i=1}^{m}\|g_iw_i^{q/s}\|_{L^{s_i}}^{s/q}\right)
\left(\prod_{i=1}^{m}\|f_i\|_{L^{p_i}}^{\frac{p_i\alpha}{nm}}\right),\\
\end{eqnarray*}
where in the last inequality we have used the fact that
$\vec{w^q}\in A_{\vec{S}}$. We observe now that
$\|g_iw_i^{q/s}\|_{L^{s_i}}^{s/q}=\|f_i\|_{p_i}^{p_i/q_i}=\|f_i\|_{p_i}^{1-\frac{\alpha
p_i}{nm}}$ and inequality ($\ref{general}$) follows immediately.

\medskip

The other implication is a consequence of the inequality ${\cal
M}_{\alpha}(\vec{f})\leq {\cal M}_{\alpha,B}(\vec{f})$ and the
boundedness result proved in \cite{M}.$\square$

\bigskip

\noindent {\it Proof of theorem $\ref{debilmaximal}$} : Let
$\Omega_{\lambda}=\{x\in \zR^n: {\cal
M}_{\alpha}\vec{f}(x)>\lambda^m\}$. By homogeneity we may assume
that $\lambda=1$. Let $K$ be a compact set contained in
$\Omega_{\lambda}$. Since $K$ is a compact set and using Vitali's
covering lemma we obtain a finite family of disjoint cubes
$\{Q_j\}$ for which
\begin{equation}\label{cubos}
1<|Q_j|^{\alpha/n}\prod_{i=1}^m \frac{1}{|Q_j|}\int_{Q_j}|f_i|,
\end{equation}
and $K\subset U_j3Q_j$. Notice that, by H\"{o}lder's inequality we
have that $\frac{u(Q)}{|Q|}\leq
\prod_{i=1}^m\left(\frac{w_i(Q)}{|Q|}\right)^{1/m}$. Then by
$(\ref{cubos})$ and H\"{o}lder's inequality at discrete level we
obtain that
\begin{eqnarray*}
u(K)^{m}&\leq&C \left(\sum_j
\frac{u(3Q_j)}{|3Q_j|}|Q_j|\right)^m\\
&\leq&C
\left(\sum_j\prod_{i=1}^m\left(\frac{1}{|3Q_j|}\int_{3Q_j}w_i\right)^{1/m}|Q_j|^{1/m}
\left(\frac{|Q_j|^{\alpha/(nm)}}{|Q_j|}\int_{Q_j}|f_i|\right)^{1/m}\right)^m\\
&\leq&C
\left(\sum_j\prod_{i=1}^m\left(\frac{|3Q_j|^{\alpha/(nm)}}{|3Q_j|}\int_{3Q_j}w_i\right)^{1/m}
\left(\int_{Q_j}|f_i|\right)^{1/m}\right)^m\\
&\leq&C \left(\sum_j\prod_{i=1}^m
\left(\int_{Q_j}|f_i|\, M_{\alpha/m}w_i\right)^{1/m}\right)^m\\
&\leq&C \prod_{i=1}^m
\int_{\zR^n}|f_i|\, M_{\alpha/m}w_i,\\
\end{eqnarray*}
and the proof concludes. $\square$

\bigskip

\noindent {\it Proof of theorem $\ref{debildebil}$}: Let $p>1$ to
be chosen later. Thus, since $L^{p,\infty}$ and $L^{p',1}$ are
associate spaces, we have that
\begin{equation*}
\|{\cal I}_{\alpha}\vec{f}\|_{L^{1/m,\infty}(u)}^{1/(pm)}=\|({\cal
I}_{\alpha}\vec{f})^{1/(pm)}\|_{L^{p,\infty}(u)}=\sup_{\|g\|_{L^{p',1}(u)}\leq
1}\int_{\zR^n}({\cal I}_{\alpha}\vec{f})^{1/(pm)} g \, u.
\end{equation*}
By theorem $\ref{p<1}$ we obtain that
\begin{equation*}
\int_{\zR^n}({\cal I}_{\alpha}\vec{f})^{1/(pm)} g \, u\leq
\int_{\zR^n}({\cal M}_{\alpha}\vec{f})^{1/(pm)} M(gu) =
\int_{\zR^n}({\cal M}_{\alpha}\vec{f})^{1/(pm)}
\frac{M(gu)}{M_{L(\log L)^{\delta}}(u)}M_{L(\log L)^{\delta}}(u),
\end{equation*}
for $\delta>0$.

\smallskip

\noindent By applying H\"{o}lder's inequality in Lorentz spaces we
obtain that
\begin{eqnarray*}
\int_{\zR^n}({\cal I}_{\alpha}\vec{f})^{1/(pm)} g \,
u&\leq&\|({\cal
M}_{\alpha}\vec{f})^{1/(pm)}\|_{L^{p,\infty}(M_{L(\log
L)^{\delta}}(u))}\left\|\frac{M(gu)}{M_{L(\log
L)^{\delta}}(u)}\right\|_{L^{p',1}(M_{L(\log L)^{\delta}}(u))}\\
\end{eqnarray*}

Now we proceed as in the linear case (see \cite{CPSS}) by taking
$p=1+\delta-2\epsilon$ with $0<2\epsilon<\delta$ which allows us
to obtain that
\begin{equation*}
\left\|\frac{M(gu)}{ M_{L(\log L)^{\delta}}(u)}
\right\|_{L^{p',1}(M(\log L)^{\delta}(u))} \leq
C\|g\|_{L^{p',1}(u)}
\end{equation*}
and taking supremum over $\|g\|_{L^{p',1}(u)}\leq 1$. $\square$

\bigskip

\noindent {\it Proof of lemma $\ref{dospesos}$}: From inequality
$(\ref{discreta})$ with $g$ replaced by $v$ and the $RH_{\infty}$
condition on $v$ we obtain that
\begin{eqnarray*}
&&\int_{\zR^n}{\cal I}_{\alpha}\vec{f}(x) u(x) v(x) dx\\
&&\qquad \leq
C\sum_{k,j}|Q_{k,j}|^{\alpha/n}\left(\frac{1}{|Q_{k,j}|}\int_{Q_{k,j}}u
v\right)\left(\prod_{i=1}^m\left(\frac{1}{|3Q_{k,j}|}\int_{3Q_{k,j}}
f_i\right)\right) |E_{k,j}|\\
&&\qquad \leq C\sum_{k,j}|Q_{k,j}|^{\alpha/n}\int_{Q_{k,j}}u
\left(\prod_{i=1}^m\left(\frac{1}{|3Q_{k,j}|}\int_{3Q_{k,j}}
f_i\right)\right)\sup_{Q_{k,j}}
v\\
&&\qquad \leq
C\sum_{k,j}|Q_{k,j}|^{\alpha/n}\left(\frac{1}{|Q_{k,j}|}\int_{Q_{k,j}}u\right)
\left(\prod_{i=1}^m\left(\frac{1}{|3Q_{k,j}|}\int_{3Q_{k,j}}
f_i\right)\right)v(Q_{k,j}). \\
\end{eqnarray*}
Since $v\in A_{\infty}$ and by the properties of the sets
$E_{k,j}$ we obtain that
\begin{eqnarray*}
&&\int_{\zR^n}{\cal I}_{\alpha}\vec{f}(x) u(x) v(x) dx\\
&&\qquad \leq
C\sum_{k,j}|3Q_{k,j}|^{\alpha/n}\left(\frac{1}{|Q_{k,j}|}\int_{Q_{k,j}}u\right)
\left(\prod_{i=1}^m\left(\frac{1}{|3Q_{k,j}|}\int_{3Q_{k,j}}
f_i\right)\right)v(E_{k,j})\\
&&\qquad \leq C\sum_{k,j}\int_{E_{k,j}}{\cal
M}_{\alpha}\vec{f}(x)Mu(x)\,v(x)\, dx\\
&&\qquad \leq C\int_{\zR^n}{\cal
M}_{\alpha}\vec{f}(x)Mu(x)\,v(x)\, dx.\ \square\\
\end{eqnarray*}

\bigskip

\noindent {\it Proof of theorem $\ref{p<1}$}: We proceed as in the
linear case (see \cite{CPSS}). We use the duality for $L^p$ spaces
for $p<1$: if $f\ge0$
\begin{equation*}
\|f\|_p=\inf\{fu^{-1}:\|u^{-1}\|_{p'}=1\}=\int fu^{-1}
\end{equation*}
for some $u\ge0$ such that $\|u^{-1}\|_{p'}=1$, with
$p'=\frac{p}{p-1}<0$. This follows from the following reverse
H\"{o}lder's inequality, which is a consequence of the
H\"{o}lder's inequality,
\begin{equation}\label{rh}
\int fg\ge \|f\|_p \|g\|_{p'}.
\end{equation}

We choose a nonnegative function $g$ such that
$\|g^{-1}\|_{L^{p'}(Mu)}=1$, and such that
\begin{equation*}
\|{\cal M}_{\alpha}\vec{f}\|_{L^{p}(Mu)}=\int {\cal
M}_{\alpha}\vec{f}\ \frac{Mu}{g}.
\end{equation*}
Let $\delta >0$. By Lebesgue differentiation theorem we get
\begin{equation*}
\|{\cal M}_{\alpha}\vec{f}\|_{L^{p}(Mu)}\ge \int {\cal
M}_{\alpha}\vec{f}\ \frac{Mu}{M_{\delta}(g)},
\end{equation*}
where $M_{\delta}(g)=M(g^{\delta})^{1/\delta}$. Then applying
lemmas $\ref{dospesos}$ and $\ref{CN}$ to the weight
$M_{\delta}(g)^{-1}$ and the reverse H\"{o}lder's inequality
$(\ref{rh})$, we obtain that
\begin{equation*}
\|{\cal M}_{\alpha}\vec{f}\|_{L^{p}(Mu)}\ge\int {\cal
I}_{\alpha}\vec{f}\ \frac{u}{M_{\delta}(g)}\ge \|{\cal
I}_{\alpha}\vec{f}\|_{L^p(u)}\|M_{\delta}(g)^{-1}\|_{L^{p'}(u)},
\end{equation*}
and everything is reduced to proved
\begin{equation*}
\|M_{\delta}(g)^{-1}\|_{L^{p'}(u)}\ge\|g^{-1}\|_{L^{p'}(Mu)}=1.
\end{equation*}
Now, the proof follows as in the linear case (see \cite{CPSS}).
Since $p'<0$, this is equivalent to prove that
\begin{equation*}
\int_{\zR^n}M_{\delta}(g)^{-p'}(x) u(x)\, dx\leq C\int
g^{-p'}(x)Mu(x)\, dx.
\end{equation*}
By choosing $\delta$ such that $0<\delta<\frac{p}{1-p}$, we have
that $-p'/\delta>1$ and the above inequality follows from the
classical weighted norm inequality of Fefferman-Stein (see
\cite{FS}). $\square$

\bigskip

\noindent {\it Proof of theorem $\ref{Welland}$}: Let $s$ be a
positive number. We split ${\cal I}_{\alpha}$ as follows
\begin{eqnarray*}
|{\cal
I}_{\alpha}\vec{f}(x)|&\leq&\int_{\sum_{i=1}^m|x-y_i|<s}\frac{\prod_{i=1}^m
|f_i(y_i)|}{(\sum_{i=1}^m|x-y_i|)^{mn-\alpha}}\,
d\vec{y}+\int_{\sum_{i=1}^m|x-y_i|\ge s}\frac{\prod_{i=1}^m
|f_i(y_i)|}{(\sum_{i=1}^m|x-y_i|)^{mn-\alpha}}\, d\vec{y}\\
&=& I_1+I_2.
\end{eqnarray*}
Let us first estimate $I_1$. Thus, if $Q_k$ is a cube centered at
$x$ with side length $2^{-k}s$, $k\in \zN\cup \{0\}$, we obtain
\begin{eqnarray*}
I_1&=&\sum_{k=0}^{\infty}\int_{2^{-k-1}s<\sum_{i=1}^m|x-y_i|\leq2^{-k}s}\frac{\prod_{i=1}^m
|f_i(y_i)|}{(\sum_{i=1}^m|x-y_i|)^{mn-\alpha}}\, d\vec{y}\\
&\leq&C\sum_{k=0}^{\infty}\frac{1}{(2^{-k}s)^{mn-\alpha}}\int_{\sum_{i=1}^m|x-y_i|\leq2^{-k}s}\left(\prod_{i=1}^m
|f_i(y_i)|\right)\, d\vec{y}\\
&\leq&C\sum_{k=0}^{\infty}\frac{1}{(2^{-k}s)^{-\alpha}}\prod_{i=1}^m\frac{1}{|Q_k|}\int_{Q_{k}}
|f_i(y_i)|\, d{y_i}\\
&\leq&C\sum_{k=0}^{\infty}\frac{1}{(2^{-k}s)^{-\alpha}}\frac{(2^{-k}s)^{\alpha-\epsilon}}
{(2^{-k}s)^{\alpha-\epsilon}}\prod_{i=1}^m\frac{1}{|Q_k|}\int_{Q_{k}}
|f_i(y_i)|\, d{y_i}\\
&\leq&C s^{\epsilon}{\cal
M}_{\alpha-\epsilon}\vec{f}(x).\\
\end{eqnarray*}
Now, we proceed to estimate $I_2$. Let $P_k$ be the cube centered
at $x$ with side length $2^k s$. Then we obtain
\begin{eqnarray*}
I_2&=&\sum_{k=0}^{\infty}\int_{2^{k}s<\sum_{i=1}^m|x-y_i|\leq2^{k+1}s}\frac{\prod_{i=1}^m
|f_i(y_i)|}{(\sum_{i=1}^m|x-y_i|)^{mn-\alpha}}\, d\vec{y}\\
&\leq&C\sum_{k=0}^{\infty}\frac{1}{(2^{k}s)^{mn-\alpha}}\int_{\sum_{i=1}^m|x-y_i|\leq2^{k+1}s}\left(\prod_{i=1}^m
|f_i(y_i)|\right)\, d\vec{y}\\
&\leq&C\sum_{k=0}^{\infty}\frac{1}{(2^{k}s)^{-\alpha}}\prod_{i=1}^m\frac{1}{|P_{k+1}|}\int_{P_{k+1}}
|f_i(y_i)|\, d{y_i}\\
&\leq&C\sum_{k=0}^{\infty}\frac{1}{(2^{k}s)^{-\alpha}}\frac{(2^{k}s)^{\alpha+\epsilon}}
{(2^{k}s)^{\alpha+\epsilon}}\prod_{i=1}^m\frac{1}{|P_{k+1}|}\int_{P_{k+1}}
|f_i(y_i)|\, d{y_i}\\
&\leq&C \frac{1}{s^{\epsilon}}{\cal
M}_{\alpha+\epsilon}\vec{f}(x).\\
\end{eqnarray*}
Collecting both estimates we obtain
\begin{equation*}
{\cal I}_{\alpha}\vec{f}(x)\leq C \left( s^{\epsilon}{\cal
M}_{\alpha-\epsilon}\vec{f}(x)+{s^{-\epsilon}}{\cal
M}_{\alpha+\epsilon}\vec{f}(x)\right),
\end{equation*}
for any $s>0$. Then, to complete the proof, we just have to
minimize the expression above in the variable $s$. $\square$

\section{Banach function spaces}\label{section5}

We introduce now some basic facts about the theory of Banach
function spaces. For more information about these spaces we refer
the reader to \cite{BS}.

Let $X$ be a Banach function space over $\zR^n$ with respect to
the Lebesgue measure. $X$ has an associate Banach function space
$X'$ for which the generalized H\"{o}lder inequality,
\begin{equation*}
\int_{\zR^n}|f(x)g(x)|\, dx\leq \|f\|_X\|g\|_{X'},
\end{equation*}
holds. Examples of Banach functions spaces are given by the
Lebesgue $L^p$ spaces, Lorentz spaces and Orlicz spaces. The
Orlicz spaces are one of the most relevant Banach function spaces,
and a brief description was given in section $\S\ref{intro}$.

Given any measurable function $f\in X$ and a cube $Q\subset
\zR^n$, we define the $X$ average of $f$ over $Q$ to be
\begin{equation*}
\|f\|_{X,Q}=\|\delta_{l(Q)}(f\chi_Q)\|_X,
\end{equation*}
where $\delta_{a}f(x)=f(ax)$ for $a>0$ and $\chi_A$ denotes the
characteristic function of the set $A$. In particular, when
$X=L^r$, $r\ge 1$, we have that
\begin{equation*}
\|f\|_{X,Q}=\left(\frac{1}{|Q|}\int_Q|f(y)|^r\right)^{1/r},
\end{equation*}
and if $X=L^B$, the Orlicz space associated to a Young function
$B$,  then
\begin{equation*}
\|f\|_{X,Q}=\|f\|_{B,Q}.
\end{equation*}

For a given Banach function space $X$, we associate the following
maximal operator defined for each locally integrable function $f$
by
\begin{equation*}
M_Xf(x)=\sup_{Q\ni x}\|f\|_{X,Q}.
\end{equation*}
If $Y_1,\dots,Y_m$ are Banach function spaces, the multilinear
version of the maximal function above is given by
\begin{equation*}
{\cal M}_{\vec{Y}}\vec{f}(x)=\sup_{Q\ni
x}\prod_{i=1}^{m}\|f_i\|_{Y_i,Q}.
\end{equation*}
Let $1<p_1,\dots,p_m<\infty$ and suppose that
$M_{Y_i}:L^{p_i}\rightarrow L^{p_i}$. From the fact that ${\cal
M}_{\vec{Y}}\vec{f}(x)\leq \prod_{i=1}^{m}M_{Y_I}f_i(x)$ and
applying H\"{o}lder's inequality we obtain that
\begin{equation*}
{\cal M}_{\vec{Y}}:L^{p_1}(\zR^n)\times\dots\times
L^{p_m}(\zR^n)\rightarrow L^p(\zR^n).
\end{equation*}

\bigskip

We define now the multilinear maximal operator associate to
certain function $\varphi$ that generalizes the multilinear
fractional maximal operator ${\cal M}_{\alpha}$. We shall assume
that the function $\varphi:(0,\infty)\rightarrow(0,\infty)$ is
essentially nondecreasing, that is, there exists a positive
constant $\rho$ such that, if $t\leq s$ then
$\varphi(t)\leq\rho\varphi(s)$. We shall also suppose that
$\lim_{t\rightarrow \infty}\frac{\varphi(t)}{t}=0$. The linear
case of the operator below was study in \cite{P3}.

\begin{definition}\label{fi}
Let $\vec{f}=(f_1,\dots, f_m)$. The multilinear maximal operator
${\mathcal M}_{\varphi}$ associated to the function $\varphi$ is
defined by
\begin{equation*}
{\cal M}_{\varphi}\vec{f}(x)=\sup_{Q\ni
x}\varphi(|Q|)\prod_{i=1}^{m}\frac{1}{|Q|}\int_Q f_i.
\end{equation*}
When $m=1$ we simply write ${\cal M}_{\varphi}=M_{\varphi}$.
\end{definition}

\bigskip

The following result is a generalized version of theorem
\ref{pesosorlicz} when $B(t)=t$. The case $m=1$ was proved in
\cite{P3}.

\medskip

\begin{theorem}\label{banach} Let $1/m<p\leq q<\infty$, $1<p_i<\infty$,
$i=1,\dots,m$, $1/p=\sum_{i=1}^{m}1/p_i$. Let $\varphi$ be a
function as in definition ($\ref{fi}$). Let $Y_i$, $i=1,\dots,m$,
be $m$ Banach function spaces such that ${\cal M}_{Y'}:
L^{p_i}\rightarrow L^{p_i}$. Suppose that $\nu, w_1,\dots,w_m$ are
weights such that, for some positive constant $C$ and for every
cube $Q$
\begin{equation}
\varphi(|Q|)|Q|^{1/q-1/p}\left(\frac{1}{|Q|}\int_Q\nu^q\right)^{1/q}\prod_{i=1}^m\|w_i^{-1}\|_{Y_i,Q}\leq
C.
\end{equation}
Then
\begin{equation}
\|{\cal M}_{\varphi}\vec{f}\nu\|_q\leq
C\prod_{i=1}^m\|f_iw_i\|_{p_i}
\end{equation}
holds for every $\vec{f}\in L^{p_1}(w_1^{p_1})\times\dots \times
L^{p_m}(w_m^{p_m})$.
\end{theorem}

When $\varphi(t)=t^{\alpha/n}$ and $Y_i=L_{A_i}$, $i=1,\dots, m$
are the Orlicz spaces associated to the Young functions $A_i$,
then we obtain theorem $\ref{pesosorlicz}$ for the case $B(t)=t$.

The proof of theorem $\ref{banach}$ follows similar arguments to
those in the proof of theorem $\ref{pesosorlicz}$. The main tools
used are an analogue of the Calder\'on-Zygmund decomposition for
${\cal M}_{\varphi}^d$ adapted to the essentially nondecreasing
function $\varphi$, the generalized H\"{o}lder's inequality and
the boundedness of ${\cal M}_{Y'}$ in the right places.

\bigskip

\begin{corollary}
Let $1/m<p<\infty$, $1<p_i<\infty$, $i=1,\dots,m$,
$1/p=\sum_{i=1}^{m}1/p_i$. Let $\varphi$ be a function as in
definition ($\ref{fi}$). Then

\noindent {\rm (i)} There exists a positive constant $C$ such
that, for every $\vec{f}=(f_1,\dots,f_m)$, and every positive
functions $u_i$
\begin{equation*}
\left(\int_{\zR^n}{\cal M}_{\varphi}\vec{f}(y)^p
(\Pi_{i=1}^mu_i(y)^{1/p_i})^p\,dy\right)^{1/p}\leq C
\prod_{i=1}^{m}\left(\int_{\zR^n}|f_i(y)|^{p_i}
M_{\varphi^p}(u_i)\right)^{1/p_i}
\end{equation*}
\noindent {\rm (ii)} If $s_i>p_i'-1$, there exists a positive
constant $C$ such that, for every $\vec{f}=(f_1,\dots,f_m)$, and
every positive functions $u_i$
\begin{equation*}
\left(\int_{\zR^n}{\cal M}_{\varphi}\vec{f}(y)^p \frac{dy} {
\left(\Pi_{i=1}^m M_{\varphi^{ps_i}}(u_i^{s_i})(y)^{1/(p_i
s_i)}\right)^p }\right)^{1/p} \leq C
\prod_{i=1}^{m}\left(\int_{\zR^n}|f_i(y)|^{p_i}\frac{dy}{u_i(y)}
\right)^{1/p_i}.
\end{equation*}
\end{corollary}
The proof of {\rm (i)} follows by applying theorem $\ref{banach}$
to the weights $\nu=\Pi_{i=1}^mu_i^{1/p_i}$,
$w_i=M_{\varphi^p}(u_i)^{1/p_i}$ and $Y_i=L^{p_i'r}$,
$1<r<\infty$.

To prove {\rm (ii)} we apply theorem $\ref{banach}$ to the weights
$\nu=\Pi_{i=1}^m M_{\varphi^{ps_i}}(u_i^{s_i})(y)^{1/(p_i s_i)}$,
$w_i=u_i^{-1/p_i}$ and $Y_i=L^{p_i'r_i}$, $r_i=(p_i-1)s_i$.


\begin{thebibliography}{AAAAAA}
\bibitem[BS]{BS} Bennett, C. and Sharpley, R.: \emph{Interpolation of operators}, Academic Press, New York
(1988).
\bibitem[BHP]{BHP}Bernardis, A., Hartzstein, S. and
Pradolini, G.: \emph{Weighted inequalities for commutators of
fractional integrals on spaces of homogeneous type}, J. Math.
Anal. Appl. 322 (2006), 825--846.
\bibitem[CCUF]{CCUF} Capone, C., Cruz Uribe, D. and  Fiorenza, A.: \emph{The fractional maximal
operator on variable $L^{p}$ spaces} preprint available in
http://www.na.iac.cnr.it.
\bibitem[C]{C} Coifman, R.:
\emph{Distribution function inequalities for singular integrals},
Proc. Nat. Acad. Sci. U.S.A. 69, No. 10 (1972), 2838--2839.
\bibitem[CF]{CF} Coifman, R. R. and Fefferman, C.: \emph{Weighted norm inequalities for
maximal functions and singular integrals},
 Studia Math. 51 (1974), 241--250.
\bibitem[CN]{CN} Cruz Uribe, D. and Neugebauer, C. J.: \emph{The structure of the reverse H\"{o}lder classes},
 Trans. Amer. Math. Soc. 347 (1995), 2941--2960.
\bibitem[CPSS]{CPSS} Carro, M. J., P\'erez, C., Soria, F. and
Soria, J.: \emph{Maximal functions and the control of weighted
inequalities for the fractional integral operator}, Indiana Univ.
Math. J.  54 (2005),  627-644.
\bibitem[CUF]{CUF} Cruz Uribe, D. and Fiorenza, A.: \emph{Endpoint estimates and weighted norm inequalities
for commutators of fractional integrals}, Publ. Mat. 47 (2003),
No. 1, 103--131.
\bibitem[CUP]{CUP} Cruz Uribe, D. and P\'erez, C.: \emph{Sharp two-weight, weak-type norm
inequalities for singular integral operators}, Math. Res. Let. 6
(1999), 417--428.
\bibitem[FS]{FS} Fefferman, C. and Stein, E. M.: \emph{Some maximal inequalities}, Amer. J. Math.
 93 (1971),107--115.
\bibitem[GCM]{GCM} Garc\'\i a-Cuerva, J. and Martell, J. M.: \emph{Two-weighted norm inequalities for maximal operators and
fractional integrals on non-homogeneous spaces}, , Indiana Univ.
Math. J. 50, No. 3, (2001), 1241-1280.
\bibitem[GCRF]{GCRF} Garc\'\i a-Cuerva, J. and Rubio de Francia, J. L.: \emph{Weighted norm inequalities and related
topics}, North Holland, Amsterdam 1985.
\bibitem[GPS]{GPS} Gorosito, O., Pradolini, G. and Salinas, O.:
\emph{Boundedness of fractional operators in weighted variable
exponent spaces with non doubling measures}, preprint.
\bibitem[G]{G} Grafakos, L.: \emph{On multilinear fractional integrals}, Studia Math. 102 (1992), 49-56.
\bibitem[GK]{GK} Grafakos, L. and Kalton, N.: \emph{Some remarks on
multilinear maps and interpolation}, Math. Ann. 319 (2001), No. 1
151-180.
\bibitem [GT]{GT} Grafakos, L. and Torres, R. H.: \emph{On multilinear singular integrals of
Calder\'on-Zygmund type}, Proceedings of the 6th international
conference on Harmonic Analysis and Partial Differential Equations
El Escorial, (2000). Publ. Mat. 2002, Vol. extra, 57--91.
\bibitem[H]{H} Hedberg, L. I.: \emph{On certain convolution
inequalities}, Porc. Amer. Math. Soc. 36 (1972), 505--510.
\bibitem[KR]{KR} Krasnosel'ski\v{i}, M. A. and Ruticki\v{i}, J.
B.:\emph{Convex function and Orlicz spaces}, Noordhoff,
Gr\"{o}ningen, 1961.
\bibitem[KS]{KS} Kenig, C. and Stein, E.: \emph{Multilinear estimates and fractional
integration}, Math. Res. Lett., 6 (1999), 1-15.
\bibitem[LOPTT]{LOPTT} Lerner, A., Ombrosi, S., P\'erez, C.,
Torres, R. and  Trujillo-Gonz\'alez, R.: \emph{New maximal
functions and multiple weights for the multilinear
Calder\'on-Zygmund theory} Adv. Math., 220 (2009) no. 4,
1222-1264.
\bibitem[M]{M} Moen, K.:  \emph{Weighted inequalities for multilinear fractional integral operators}, to appear in Collectanea Mathematica
\bibitem [MW]{MW} Muckenhoupt, B. and Wheeden, R.: \emph{Weighted norm inequalities for fractional
integrals}, Trans. Amer. Math. Soc., 192 (1974), 261-274.
\bibitem[O]{O} O'Neil, R.: \emph{Fractional integration in Orlicz
spaces}, Trans. Amer. Math. Soc. 115 (1963) 300, 328.
\bibitem[P1]{P1} P\'erez, C.: \emph{On sufficient conditions for the
boundedness of the Hardy-Littlewood maximal operator between
weighted $L^p$-spaces with different weights }, Proc. London Math.
Soc. (3) 71 (1995) 135, 157.
\bibitem[P2]{P2} P\'erez, C.: \emph{Sharp $L^p$-weighted Sobolev inequalities},
Ann. Inst. Fourier 45 (3) (1995) 809, 824.
\bibitem[P3]{P3} P\'erez, C.: \emph{Two weighted inequalities for potential and fractional type maximal operators},
Indiana Univ. Math. 43 (2) (1994) 663, 683.
\bibitem[P4]{P4} P\'erez, C.: \emph{Endpoint estimates for commutators of singular integral
operators}, J. Funct. Anal. 127 (1995) 163--185.
\bibitem[P5]{P5} P\'erez, C.: \emph{Sharp estimates for commutators of singular integrals
via iterations of the Hardy–Littlewood maximal function}, J.
Fourier Anal. Appl. 3 (6) (1997) 743–756.
\bibitem[PT]{PT} P\'erez, C. and Trujillo-Gonz\'alez, R.: \emph{Sharp weighted estimates
for multilinear commutators}, J. London Math. Soc., 65 (2), (2002)
672--692 .
\bibitem[RR]{RR} Rao, M. M. and Ren, Z. D.:\emph{ Theory of Orlicz
spaces}, Monogr. Textbooks Pure Appl. Math., 146, Dekker, New
York, 1991.
\bibitem[W]{W} Welland, G. V.: \emph{Weighted norm inequalities for fractional
integrals}, Proc. Amer. Math. Soc. 51 (1975), 143--148.
\end{thebibliography}
\end{document}